# Generalizing Wallis' formula

D Huylebrouck

**Abstract.** The present note generalizes Wallis' formula, $\frac{\pi}{2} = \frac{2}{1} \cdot \frac{2}{3} \cdot \frac{4}{3} \cdot \frac{4}{5} \cdot \frac{6}{5} \cdot \frac{6}{7} \cdot \ldots$, using the Euler-Mascheroni constant g and the Glaisher-Kinkelin constant A:

$$\frac{4^{\gamma}}{2^{\ln 2}} = \frac{\sqrt{2}}{1} \cdot \frac{\sqrt{2}}{\sqrt[3]{3}} \cdot \frac{\sqrt[4]{4}}{\sqrt[3]{3}} \cdot \frac{\sqrt[4]{4}}{\sqrt[5]{5}} \cdot \frac{\sqrt[6]{6}}{\sqrt[5]{5}} \cdot \frac{\sqrt[6]{6}}{\sqrt[7]{7}} \cdot \ldots \quad \text{and} \quad \left(\frac{2^2 \pi e^{\gamma}}{A^{12}}\right)^{\frac{\pi^2}{6}} = \frac{\sqrt[4]{2}}{1} \cdot \frac{\sqrt[4]{2}}{\sqrt[9]{3}} \cdot \frac{\sqrt[16]{4}}{\sqrt[9]{3}} \cdot \frac{\sqrt[16]{4}}{\sqrt[25]{5}} \cdot \frac{\sqrt[36]{6}}{\sqrt[25]{5}} \cdot \frac{\sqrt[36]{6}}{\sqrt[49]{7}} \cdot \ldots$$

Wallis' formula, named after the English mathematician John Wallis (1616 –1703), is popular in many calculus courses (see [**2**], [**1**] p. 338]). It is a slowly convergent product, but its importance is historic and aesthetic. The present paper proposes two similar equally pleasing formulas, in a rather straighforward way, without using the gamma function for generating these product formulas. Perhaps some readers will take up the challenge of finding even easier proofs on the level of a calculus course, similar to those for Wallis' formula.

The Dirichlet eta function is defined for any complex number s with real part > 0 by:

$$\eta(s) = \frac{1}{1^s} - \frac{1}{2^s} + \frac{1}{3^s} - \frac{1}{4^s} + \frac{1}{5^s} + \ldots \; .$$

For $s = 0$ it lead Y. L. Yung and J. Sondow to a remarkably elegant proof for Wallis' formula $\frac{\pi}{2} = \frac{2}{1} \cdot \frac{2}{3} \cdot \frac{4}{3} \cdot \frac{4}{5} \cdot \frac{6}{5} \cdot \frac{6}{7} \cdot \ldots$ (see [**2**]), and their proof will be adapted here to other values of s.

**Theorem.** *For appropriate values of s (and if $\sqrt[n^0]{2n}$ is interpreted as 2n),*

$$e^{2\eta'(s)} = \frac{\sqrt[2^s]{2}}{\sqrt[1^s]{1}} \cdot \frac{\sqrt[2^s]{2}}{\sqrt[3^s]{3}} \cdot \frac{\sqrt[4^s]{4}}{\sqrt[3^s]{3}} \cdot \frac{\sqrt[4^s]{4}}{\sqrt[5^s]{5}} \cdot \frac{\sqrt[6^s]{6}}{\sqrt[5^s]{5}} \cdot \frac{\sqrt[6^s]{6}}{\sqrt[7^s]{7}} \cdot \ldots \qquad (*)$$

*Proof.* By definition,

$$\eta(s) = \frac{1}{1^s} - \frac{1}{2^s} + \frac{1}{3^s} - \frac{1}{4^s} + \frac{1}{5^s} + \ldots \qquad \text{(for Re}(s) > 0)$$

$$= \frac{1}{2} + \frac{1}{2}\left[\left(\frac{1}{1^s} - \frac{1}{2^s}\right) + \left(-\frac{1}{2^s} + \frac{1}{3^s}\right) + \left(\frac{1}{3^s} - \frac{1}{4^s}\right) + \left(-\frac{1}{4^s} + \frac{1}{5^s}\right) + \ldots\right] \qquad \text{(for Re}(s) > -1)$$

Thus: $\eta'(s) = \frac{1}{2}\left[\left(-\frac{\ln(1)}{1^s} + \frac{\ln(2)}{2^s}\right) + \left(\frac{\ln(2)}{2^s} - \frac{\ln(3)}{3^s}\right) + \left(-\frac{\ln(3)}{3^s} + \frac{\ln(4)}{4^s}\right) + \left(\frac{\ln(4)}{4^s} - \frac{\ln(5)}{5^s}\right) + \ldots\right]$

$$= \frac{1}{2}\ln\left(\frac{2^{1/2^s}}{1^{1/1^s}} \cdot \frac{2^{1/2^s}}{3^{1/3^s}} \cdot \frac{4^{1/4^s}}{3^{1/3^s}} \cdot \frac{4^{1/4^s}}{5^{1/5^s}} \ldots\right),$$

and so $e^{2\eta'(s)} = \frac{2^{1/2^s}}{1^{1/1^s}} \cdot \frac{2^{1/2^s}}{3^{1/3^s}} \cdot \frac{4^{1/4^s}}{3^{1/3^s}} \cdot \frac{4^{1/4^s}}{5^{1/5^s}} \ldots$ ∎

In the case $s = 0$, the litterature tells us $\eta'(0) = \frac{1}{2}\ln\left(\frac{\pi}{2}\right)$, and this leads to the familiar Wallis formula $\frac{\pi}{2} = \frac{2}{1}\cdot\frac{2}{3}\cdot\frac{4}{3}\cdot\frac{4}{5}\cdot\frac{6}{5}\cdot\frac{6}{7}\cdot\ldots$.

In the case $s = 1$, $\eta'(1)$ can be obtained from the property $\sum_{k=1}^{\infty}(-1)^k\frac{\ln(k)+\gamma}{k} = -\frac{1}{2}(\ln 2)^2$ (see [5]). Indeed, using the Euler-Mascheroni constant γ, it follows that

$$\eta'(1) = \sum_{k=1}^{\infty}(-1)^k\frac{\ln(k)}{k} = -\gamma\sum_{k=1}^{\infty}(-1)^k\frac{1}{k} - \frac{1}{2}(\ln 2)^2$$
$$= -\gamma(-\ln 2) - \frac{1}{2}(\ln 2)^2 = \frac{1}{2}(2\gamma - \ln 2)\ln 2 = \frac{1}{2}\ln 2^{(2\gamma - \ln 2)},$$

Substitution in (*) leads to:

$$2^{(2\gamma - \ln 2)} = \frac{2^{1/2}}{1^{1/1}}\cdot\frac{2^{1/2}}{3^{1/3}}\cdot\frac{4^{1/4}}{3^{1/3}}\cdot\frac{4^{1/4}}{5^{1/5}}\cdots,$$

and this is the formula given in the abstract.

In the case $s = 2$, we need the Glaisher-Kinkelin constant $A$, given by

$$A = \lim_{n\to+\infty}\frac{H(n)}{n^{\frac{n^2}{2}+\frac{n}{2}+\frac{1}{12}}\cdot e^{-\frac{n^2}{4}}} = 1.28242\ldots$$

where the 'hyperfactorial' $H(n)$ is defined by $H(n) = \prod_{k=1}^{n} k^k$.

Now combining the expression for $\eta'(s)$ (see [3]) for $s = 2$ and a well-known result using $A$ (see [4]), implies that

$$\eta'(2) = \frac{1}{2}\frac{\pi^2}{6}\ln 2 + \frac{1}{2}\frac{\pi^2}{6}[\gamma + \ln(2\pi) - 12\ln A] = \frac{1}{2}\frac{\pi^2}{6}\ln\left(\frac{2^2\pi\,e^\gamma}{A^{12}}\right).$$

Substitution in (*) leads to:

$$\left(\frac{2^2\pi\,e^\gamma}{A^{12}}\right)^{\frac{\pi^2}{6}} = \frac{2^{1/2^2}}{1^{1/1^2}}\cdot\frac{2^{1/2^2}}{3^{1/3^2}}\cdot\frac{4^{1/4^2}}{3^{1/3^2}}\cdot\frac{4^{1/4^2}}{5^{1/5^2}}\cdots$$

and this is the formula given in the abstract.

**References**
1. Jonathan M. Borwein and Peter B. Borwein, *Pi and the AGM: A Study in Analytic Number Theory and Computational Complexity*, ISBN: 0-471-31515-X, 432 pages, July 1998.
2. J. Sondow and E. W. Weisstein, Wallis' Formula, *MathWorld - A Wolfram Web Resource*, http://mathworld.wolfram.com/WallisFormula.html.
3. E. W. Weisstein, Dirichlet Eta Function, *MathWorld - A Wolfram Web Resource*, http://mathworld.wolfram.com/DirichletEtaFunction.html.
4. R. G. Wilson, Decimal Expansion of Zeta'(2) (the first derivative of the zeta function at 2), Sequence A000129, *The On-line Encyclopedia of Integer Sequences*, http://oeis.org/A073002.
5. Equality with Euler–Mascheroni constant:
http://math.stackexchange.com/questions/100871/equalitywith-eulermascheroni-constant.